\documentclass{article}

\usepackage{amsmath,amssymb,amsthm,setspace,dsfont,hyperref,kpfonts}
\usepackage[T1]{fontenc}

\hypersetup{
   colorlinks,
    citecolor=black,
    filecolor=black,
    linkcolor=black,
    urlcolor=black
    }

\renewcommand{\thispagestyle}[1]{} 

\usepackage{fancyhdr}
\usepackage{lastpage}

\newtheorem{theorem}{Theorem}[section]
\newtheorem{lemma}[theorem]{Lemma}

\newtheorem{proposition}[theorem]{Proposition}
\newtheorem{corollary}[theorem]{Corollary}

\theoremstyle{definition}

\numberwithin{equation}{section}

\theoremstyle{remark}

\newtheorem*{note*}{Note}

\makeatletter
\newcommand{\ls}{\leqslant}
\newcommand{\gr}{\geqslant}
\newcommand{\prend}{$\hfill \Box$}
\makeatother

\begin{document}

\small

\title{Neighborhoods on the Grassmannian of marginals with bounded isotropic constant}

\author{Grigoris Paouris\thanks{Supported by the A. Sloan foundation, BSF grant 2010288 and the NSF
CAREER-1151711 grant;} \and Petros Valettas}

\date{April 19, 2014}

\maketitle

\begin{abstract}
\footnotesize We show that for any isotropic log-concave probability
measure $\mu$ on $\mathbb R^n$, for every $\varepsilon >0$, every
$1\ls k\ls \sqrt{n}$ and any $E\in G_{n,k}$ there exists $F\in
G_{n,k}$ with $d(E,F)<\varepsilon$ and $L_{\pi_F\mu}<
C/\varepsilon$.
\end{abstract}

\bigskip

\medskip

\section{Introduction}

\noindent Let $K$ be a symmetric convex body in $\mathbb R^{n}$ of volume $1$, $|K|=1$. The Hyperplane conjecture
(posed by J. Bourgain in
\cite{Bou}) claims that there exists a universal constant $c>0$ and a unit vector $\theta$ ($ \theta \in S^{n-1}$)
such that
$$ | K\cap \theta^{\perp} | \gr c, $$
where $| \cdot |$ stands for the $n$-dimensional volume. K. Ball in
\cite{Ball} showed that the question has an equivalent formulation
in the more general setting of
$\log$-concave measures. A Borel probability measure $\mu$ on
$\mathbb R^n$ is called log-concave if for any compact sets $A, B$
in $\mathbb R^n$ we have $$
  \mu((1-\lambda)A+\lambda B)\gr \mu(A)^{1-\lambda} \mu(B)^\lambda$$ for
  all $\lambda\in (0,1)$. The
measure $\mu$ is called centered if $\int \langle x,y\rangle\,
d\mu(x)=0$ for all $y\in \mathbb R^n$. The covariance matrix for a
centered measure $\mu$ is defined as: $${\rm
Cov}(\mu)_{ij}=\int_{\mathbb R^n} x_i x_j \, d\mu(x), \quad
i,j=1,\ldots,n.$$ A log-concave probability measure $\mu$ on
$\mathbb R^n$ is called isotropic if it is centered and its
covariance matrix is the identity.

\noindent The isotropic constant of a centered log-concave probability measure
on $\mathbb R^n$ is defined as:
\begin{equation}
  L_\mu:= \|\mu\|_\infty^{1/n} [\det{\rm Cov}(\mu)]^{\frac{1}{2n}},
\end{equation} where $\|\mu\|_\infty=\|f_\mu\|_\infty$ and $f_\mu$ is the density function
of $\mu$. So if $\mu$ is isotropic then
$L_\mu=\|\mu\|_\infty^{1/n}$. The Hyperplane
Conjecture can be formulated equivalently as follows: There exists a constant $C>0$ such that
for all $n$ and any $\mu$ $\log$-concave isotropic probability measure in $\mathbb
R^n$, $L_\mu<C$. A classical
reference on the problem is \cite{MiPa}. For a more detailed
exposition on recent developments see \cite{BGVV}.

The first non-trivial bound on this question was given by Bourgain
in \cite{Bou} proving that $L_{\mu}\ls c_1\sqrt[4]{n}\log n$ for all
isotropic measures on $\mathbb R^n$. In \cite{Kl} Klartag removed
the logarithmic term (see also \cite{KM} for an alternative approach
and \cite{Vri} for further refinements). Throughout this note, all
constants $c, C, C', \ldots$ denote positive, dimension-independent
numerical constants, whose value may change from line to line. We
write $A\simeq B$ to denote $c\ls A/B\ls C$ for some numerical
constants $c,C>0$.

Actually Klartag's approach developed in \cite{Kl} gave an
affirmative answer to the ``isomorphic version'' of the Hyperplane conjecture in the setting of convex
bodies:

\smallskip

\noindent {\bf Theorem} (Klartag) {\it For any $\varepsilon>0$ and
any convex body $K$ on $\mathbb
  R^n$ there exists a convex body $T$ on $\mathbb R^n$ which
  satisfies:
  \begin{itemize}
    \item $d_G(K, T)<1+\varepsilon$ and
    \item $L_T< C/\sqrt{\varepsilon}$.
  \end{itemize}}

\smallskip

In this note we show that a result of the same flavor holds true for
marginals of a log-concave isotropic probability measure. Marginals of an
isotropic log-concave measure are also log-concave from Pr\'{e}kopa-Leindler
inequality \cite{Prek} and isotropic. However, if $\mu$ has bounded
isotropic constant, it is not known whether its marginals have also
bounded isotropic constant. It's not hard to show that this is an
equivalent formulation of the problem (see \cite{Pa3}, Proposition 5.3).
Our main result states that even if for a given marginal one can not decide if the isotropic
constant is bounded, there exists another one
``close'' to it which has bounded isotropic constant. To formulate
the statement precisely let us fix the distance $d$ on $G_{n,k}$ as:
\begin{align}
  E, F\in G_{n,k} , \quad  d(E,F):=\inf \{\|I-U\|_{\rm op} : U\in O(n),
  U(E)=F\}.
\end{align} Under this notation our result reads as follows:

\begin{theorem} \label{thm:main-result}
  Let $\mu$ be isotropic log-concave probability measure on $\mathbb R^n$ and let $1\ls k\ls \sqrt{n}$.
  For every $\varepsilon >0$ and any $E\in
  G_{n,k}$ there exists $F\in G_{n,k}$ with $d(E, F)<\varepsilon$
  such that \begin{equation} \label{eq:1.4}
    L_{\pi_F\mu}\ls C/\varepsilon,
  \end{equation} where $ C>0$ is an absolute constant. Additionaly,
  if $L_\mu$ is bounded then we can take $1\ls k\ls n-1$.
\end{theorem}

Apart from the obvious similarities to Klartag's result, there are
some significant differences that we would like to point out. First, the
way of measuring distance is different.
Second, both measures we consider
here, the given and the perturbed one, are isotropic while in Klartag's result are not.
Third, the dependence on $\varepsilon$ is weaker than Klartag's result.
However, in the last section of this note we show that if one can
prove the above statement with a dependence $o(1/\varepsilon)$ then
Hyperplane Conjecture will follow. In other words, if Hyperplane
Conjecture turns out to be false, the dependence on $\varepsilon$ in
Theorem \ref{thm:main-result} is optimal.

Let us also comment on the the range of the dimension that appears
on the theorem. In order to find a marginal, close to any other
marginal, with bounded isotropic constant we have to know that there
exists at least one with this property on dimension $k$. In general,
the dimension that one knows the existence of a marginal with bounded
isotropic constant is up to order $\sqrt{n}$. On the other hand,
if $L_\mu$ is bounded, one can find marginal with bounded
isotropic constant in all lower dimensions. For a more detailed
analysis on the range of dimension that Theorem
\ref{thm:main-result} holds see Theorem \ref{thm:main-2} and the
remark after.

The paper is organized as follows: In Section 2 we give the
necessary background for the Grassmann manifold and some facts from
convexity that we are going to use. In Section 3 we briefly refer to
some standard facts about log-concave measures and we proceed in the
proof of the main result. Finally, in Section 4 we conclude with the
optimality on the theorem and some further remarks.

\section{Grassmann manifold and convexity}

\noindent {\bf \S2.1. The Grassmann manifold.} The Grassmann
manifold consisting by all $k$-dimensional subspaces $F$ of $\mathbb
R^n$ is denoted by $G_{n,k}$.
We work in $G_{n,k}$ equipped with a metric $\rho$ which
is induced by some unitarily invariant ideal norm on $L(\ell_2^n)$
(see \cite{Sz2} for details). A typical example of such a metric is
$\sigma_\infty(E_1,E_2):=\|P_{E_1}-P_{E_2}\|_{\rm op}$ induced by
the operator norm under the embedding $F\mapsto P_F$ of the
Grassmann manifold into $L(\ell_2^n)$. Another example of an
unitarily invariant metric on $G_{n,k}$ is: $d(E,F)=\inf\{\|I-U\| :
U\in O(n), \, U(E)=F\}$. Note that $\sigma_{\infty}$ and $d$ are
Lipschitz equivalent metrics, i.e. $\sigma_\infty(E,F)\ls d(E,F)\ls
\sqrt{2} \sigma_\infty(E,F)$, for all $E,F\in G_{n,k}$.
Under the orthogonal group action over $G_{n,k}$ we
get that $G_{n,k}$ becomes a homogeneous  space, thus there exists a
unique probability measure $\nu_{n,k}$ which is invariant under this
action -- the so-called Haar measure.

We will need an entropy estimate for $G_{n,k}$ with respect
to metrics described above. For a compact set $A$ in a
metric space $(X,d)$ and any $\varepsilon>0$ the
$\varepsilon$-entropy of $A$ denoted by $N(A,d, \varepsilon )$ is
the minimum number of balls of radius $\varepsilon$ required to
cover $A$. In \cite{Sz1} (see also \cite{Sz2}) Szarek proved the
following result:

\begin{theorem} \label{thm:Grass-entropy} Let $1\ls k\ls n-1$ and let $\rho$ be a metric on $G_{n,k}$ induced by some
unitarily invariant ideal norm on $L(\ell_2^n)$ normalized such that
${\rm diam}(G_{n,k},\rho)=1$. Then, for any $0<\varepsilon<1$ we
have:
  \begin{equation}
    (c_1/\varepsilon)^{k(n-k)}\ls
    N(G_{n,k}, \rho, \varepsilon)\ls (c_2/\varepsilon)^{k(n-k)},
  \end{equation} where $c_1,c_2>0$ are absolute constants.
\end{theorem}

See also \cite{Paj} for an alternative proof of this result. Using
Theorem \eqref{thm:Grass-entropy}, the invariance of the Haar
measure and the sub-additivity of the measure we can immediately
conclude the following:

\begin{corollary}\label{cor:measure-ball-grass}
  Let $1\ls k\ls n-1$, $\rho$ be a metric as above on $G_{n,k}$ and let $\delta\in (0,1)$.
  Then, for any $F\in G_{n,k}$ we have:
\begin{equation}\label{eq:2.2}
(\delta/c_2)^{k(n-k)}\ls \nu_{n,k}(B_\rho(F, \delta))\ls
(\delta/c_1)^{k(n-k)},
\end{equation} where $B_\rho(F,\delta)=\{E\in G_{n,k} :
\rho(F,E)<\delta\}$.
\end{corollary}

We are going to use only the left hand-side inequality in
\eqref{eq:2.2}.

\medskip

\noindent {\bf \S2.2. Affine and dual affine quermassintegrals.} The
next affine invariants of any body $K$ were initially introduced by
Lutwak (under different normalization) in \cite{Lut1}. The
definition we expose here follows \cite{DP}. Let $K$ be a convex
body in $\mathbb R^n$. For any $1\ls k\ls n-1$ we define the {\it
$k$-th dual affine quermassintegral} of $K$ as:
  \begin{equation}
    \tilde{\Phi}_{[k]}(K):=\left(\int_{G_{n,k}}|K\cap F^\perp|^n \,
    d\nu_{n,k}(F)\right)^{\frac{1}{kn}},
  \end{equation} where $|\cdot|$ stands for the volume.

Grinberg proved in \cite{Grin1} that these quantities are invariant
under volume preserving linear transformations, as conjectured by
Lutwak. The following inequality was proved by Busemann and Strauss \cite{BS} and independently by Grinberg in \cite{Grin1}.

\begin{theorem} \label{thm:Grinberg-ineq-dual-aff} For any $1\ls k\ls n-1$ and any convex body
$K$ in $\mathbb R^n$ of volume $1$, we have:
  \begin{equation} \int_{G_{n,k}}|K\cap E|^n \, d\nu_{n,k}(E)\ls \int_{G_{n,k}}|D_n\cap E|^n\, d\nu_{n,k}(E).
  \end{equation}
\end{theorem}

\smallskip

Here $D_n$ stands for the Euclidean ball of volume $1$. Note that
according to the above notation this inequality can be equivalently
rewritten as $\tilde \Phi_{[k]}(K)\ls \tilde \Phi_{[k]}(D_n)$ for all convex bodies of volume $1$.
We also have the following asymptotic estimate:
$$\tilde \Phi_{[k]}(D_n)=\omega_n^{1/n}
    \left(\frac{\omega_{n-k}}{\omega_n}\right)^{1/k} \simeq 1,$$
where $\omega_m$ denotes the volume of the Euclidean ball on
$\mathbb R^m$ of radius $1$.

\medskip

\noindent Next, we give the definition of the affine
quermassintegrals: Let $K$ be a convex body in $\mathbb R^n$. For
any $1\ls k\ls n-1$ we define the {\it $k$-th affine
quermassintegral} of $K$ as:
\begin{equation}\label{eq:2.5}
    \Phi_{[k]}(K):=\left( \int_{G_{n,k}} |P_FK|^{-n}\,
    d\nu_{n,k}(F) \right)^{-\frac{1}{kn}}.
  \end{equation}
\noindent We have that $K\mapsto
\Phi_{[k]}(K)$ is homogeneous of order $1$, that is if $\lambda>0$
then $\Phi_{[k]}(\lambda K)=\lambda \Phi_{[k]}(K)$.

\smallskip

These quantities were introduced by Lutwak in \cite{Lut2}
under different normalization. The definition given here follows
again \cite{DP}. Grinberg also proved in \cite{Grin1} that these
quantities are invariant under volume preserving affine
transformations. The quantity ${\Phi_{[k]}}(\cdot)$ is an affine
version of the classic quermassintegral. We refer to the book
\cite{Sch} for additional basic facts from Convex Geometry.

Lutwak conjectured that for any convex body $K$ on $\mathbb R^n$ of
volume $1$ one must have $\Phi_{[k]}(K) \gr \Phi_{[k]}(D_n)$. An
isomorphic version of this estimate was verified by P. Pivovarov and the first name author
in \cite[Theorem
5.1]{PP}. This estimate is crucial for our argument.
Therefore, we give a brief sketch of proof of this result
for reader's convenience.  The argument makes use of Theorem
\ref{thm:Grinberg-ineq-dual-aff}, Blaschke-Santal\'{o}'s inequality
\cite{Sa} and reverse Santal\'{o} inequality due to Bourgain and V.
Milman \cite{BM}.

\begin{theorem} \cite{PP} \label{thm:lower-bd-PP} Let $K$ be a centrally symmetric convex body
  on $\mathbb R^n$. Then, for all $1\ls k\ls n-1$ one has:
  \begin{equation}
    \Phi_{[k]}(K)\gr c_1|K|^{1/n}\sqrt{\frac{n}{k}},
  \end{equation} where $c_1>0$ is an absolute constant.
\end{theorem}

\noindent {\it Proof (Sketch).} Using reverse Santal\'{o} inequality
we may write:
$$\Phi_{[k]}(K)\gr c\omega_k^{2/k} \left( \int_{G_{n,k}} |K^\circ \cap F|^n \, d\nu_{n,k}(F)
\right)^{-\frac{1}{kn}}.$$ On the other hand by Theorem
\ref{thm:Grinberg-ineq-dual-aff} we have:
$$\int_{G_{n,k}}|K^\circ \cap F|^n \, d\nu_{n,k}(F)\ls
|K^\circ |^k \int_{G_{n,k}}|D_n\cap F|^n \,
d\nu_{n,k}(F)=|K^\circ|^k \omega_n^{-k}\omega_k^n.$$ Hence, by
Santal\'{o}'s inequality we obtain:
$$\left(\int_{G_{n,k}} |K^\circ \cap F|^n \, d\nu_{n,k}(F)
\right)^{\frac{1}{kn}} \ls |K^\circ |^{1/n}
\omega_n^{-1/n}\omega_k^{1/k}\ls |K|^{-1/n}\omega_n^{1/n}
\omega_k^{1/k}.$$ The result follows. \prend

\section{Proof of the main result}

We first collect some known results from the theory of $\log$-concave probability measures that we will need for the proof.
For any log-concave probability measure $\mu$ and every $q\gr 1$ we define the
{\it $L_q$-centroid body of $\mu$}, denoted by $Z_q(\mu)$ through
its support function:
\begin{equation}\label{eq:3.1}
  h_{Z_q(\mu)}(y):=\left(\int_{\mathbb R^n} |\langle x, y\rangle|^q\,
  d\mu(x)\right)^{1/q}, \quad y\in \mathbb R^n.
\end{equation}
By the definition of the marginal we have that for every subspace $F$ of $\mathbb R^{n}$,
\begin{align} \label{eq:3.6} P_F(Z_q(\mu))=Z_q(\pi_F\mu)
\end{align}
where $P_{F}$ is the orthogonal projection onto $F$ and $\pi_{F} \mu$ the marginal of $\mu$ on $F$.
Note that the isotropicity of a centered $\mu$ can be equivalently described by
the condition $Z_2(\mu)=B_2^n$. The following estimate for the volume of  $L_n$-centroid
body has been proved in  \cite[Proposition
3.7]{Pa2}:
\begin{equation}\label{eq:3.2}
  |Z_n(\mu)|^{1/n}\simeq \frac{1}{f_\mu(0)^{1/n}}.
\end{equation} Moreover, for any centered log-concave measure by Fradelizi's theorem \cite{Frad} we know
that:
\begin{equation} \label{eq:3.3}
  f_\mu(0)\ls \|f_\mu\|_\infty \ls e^n f_{\mu}(0)
\end{equation} and in view of the definition of the isotropic
constant we may write:
\begin{equation}\label{eq:3.4}
  |Z_n(\mu)|^{1/n} \simeq \frac{[\det {\rm
  Cov}(\mu)]^{\frac{1}{2n}}}{L_\mu}\ls c_1 [\det {\rm
  Cov}(\mu)]^{\frac{1}{2n}},
\end{equation} where in the last inequality we have used the fact
that $L_\mu\gr c>0$ for all probability measures $\mu$ (see
\cite{MiPa}). Moreover, for the full range of $1\ls q\ls n$ we have that \cite{Pa1}:
\begin{equation}
  |Z_q(\mu)|^{1/n}\ls c\sqrt{\frac{q}{n}} [\det {\rm Cov}(\mu)]^{\frac{1}{2n}}.
\end{equation}
For our purpose we will like to know the reverse inequality. Of course if the reverse inequality was
known for all $q\ls n$ then (just apply for $q=n$ and use \eqref{eq:3.2} and \eqref{eq:3.3}) the Hyperplane
Conjecture would follow. We introduce (for $\beta \gr 1$) the auxiliary parameter $q_v(\mu,
\beta)$ as follows:
\begin{equation}
q_v(\mu, \beta):=\max\left \{q\ls n: |Z_q(\mu)|^{1/n}\gr
\frac{1}{\beta}\sqrt{\frac{q}{n}} [\det {\rm
  Cov}(\mu)]^{\frac{1}{2n}}\right\}.
  \end{equation}
The results of Klartag and E. Milman in \cite{KM} imply that for any
$\log$-concave measure $\mu$ on $\mathbb R^n$ inequality can be reversed
(up to absolute constants) for $q\ls \sqrt{n}$ (see also
\cite{Vri} for subsequent refinements). In our notations their results give the following
\begin{theorem} \label{KM} \cite{KM} There exists an absolute constant $c>0$ such that for every centered
$\log$-concave probability measure $\mu$,
$$ q_{v}(\mu, c) \gr \sqrt{n} . $$
\end{theorem}

The Lutwak-Yang-Zhang inequalities \cite{LYZ} (see also
\cite{PP2} for the measure theoretic version we use here) say
that for all $1\ls q\ls n$,
\begin{equation}\label{eq:3.8}
  |Z_q(\mu)|^{1/n} \gr c\sqrt{\frac{q}{n}} \frac{[\det {\rm Cov}(\mu)]^{\frac{1}{2n}}}{L_\mu}.
\end{equation} In our notation \eqref{eq:3.8} implies that
\begin{equation} \label{eq:3.9-}
 q_{v}(\mu, cL_{\mu}) \gr n .
\end{equation}
Also, since $p\mapsto |Z_{p}(\mu)|^{\frac{1}{n}}$ is increasing by H\"older's inequality, we have that for every $t,\beta\gr 1$,
\begin{equation}
 q_{v} (\mu, t\beta ) \gr cq_{v} (\mu, \beta) t^{2} .
\end{equation}
For any log-concave probability measure $\mu$ on $\mathbb R^n$ we introduce the
averages:
\begin{equation}
  {\cal A}_{[k]}(\mu):=\left(\int_{G_{n,k}} f_{\pi_F\mu}(0)^n\,
  d\nu_{n,k}(F)\right)^{\frac{1}{kn}}.
\end{equation} The next Lemma shows that the above quantities are closely
related to the affine quermassintegrals and in turn to the volume of
the $L_k$-centroid body of $\mu$.

\begin{lemma} \label{lem:A_k-Phi-formula}
  Let $\mu$ be a $log$-concave isotropic measure in $\mathbb R^n$. For all $1\ls k\ls n-1$ we have:
  \begin{equation}\label{eq:3.9}
    {\cal A}_{[k]}(\mu)\simeq \left(\int_{G_{n,k}} L_{\pi_F\mu}^{kn}\, d\nu_{n,k}(F)\right)^{\frac{1}{kn}}
    \simeq \Phi_{[k]}(Z_k(\mu))^{-1}.
  \end{equation} In particular, we have that:
  \begin{equation}\label{eq:3.10}
    {\cal A}_{[k]}(\mu)\ls c_1\sqrt{\frac{k}{n}}
    \frac{1}{|Z_k(\mu)|^{1/n}},
  \end{equation} where $c_1>0$ is an absolute constant.
\end{lemma}

\noindent {\it Proof.} The first equivalence follows directly from
 \eqref{eq:3.3} applied to $\pi_F\mu$ while the second
equivalence from \eqref{eq:3.2}, \eqref{eq:3.6} and \eqref{eq:2.5}
for $K=Z_k(\pi_F\mu)$. The estimate \eqref{eq:3.10} follows from \eqref{eq:3.9} and by
Theorem \ref{thm:lower-bd-PP} applied to $Z_k(\mu)$. \prend

\medskip

\noindent So, if $\mu$ is isotropic then \eqref{eq:3.10}  and \eqref{eq:3.9-} implies that
\begin{align}\label{eq:3.13} {\cal A}_{[k]}(\mu)\ls C\beta  \  {\rm if } \ 1\ls k \ls q_{v}(\mu, \beta) \ \ \ {\rm and } \  \ \ {\cal A}_{[k]} (\mu) \ls c L_{\mu} \ {\rm if } \ k\ls n-1.\end{align}
An application of Markov's inequality yields the following
large deviation estimate:

\begin{proposition} \label{prop:large-dev-iso}
  Let $\mu$ be an isotropic log-concave probability measure on $\mathbb R^n$, $\beta\gr 1$ and
  $1\ls k\ls q_v(\mu, \beta)$. Then, we have:
  \begin{equation}
    \nu_{n,k}(\{F\in G_{n,k} : L_{\pi_F\mu} \gr C\beta t\})\ls
    t^{-kn},
  \end{equation} for all $t>1$.
\end{proposition}

Now we are ready to prove the following:

\begin{theorem}\label{thm:main-2}
  Let $\mu$ be an isotropic log-concave probability measure on $\mathbb
  R^n$ and let $\beta\gr 1$. For any $1\ls k\ls q_v(\mu, \beta)$, any $E\in G_{n,k}$ and
  every $\varepsilon \in (0,1)$, there exists $F\in G_{n,k}$ such that $\rho(E,F)<\varepsilon$ and
  $$L_{\pi_F \mu}< \frac{C \beta}{ \varepsilon^{1-\frac{k}{n}}},$$ where $C>0$ is an absolute constant and $\rho$ is
  a metric on $G_{n,k}$ as in Theorem \ref{thm:Grass-entropy}.
\end{theorem}

\noindent{\it Proof.} Let $\varepsilon \in (0,1), \, \beta\gr 1$,
$1\ls k\ls q_v(\mu, \beta)$ and let $E\in G_{n,k}$. For any
$t>1$ consider the set $A_t:=\{F\in G_{n,k} : L_{\pi_F\mu}\gr C_1
t \beta\}$. Proposition \ref{prop:large-dev-iso}
implies that $\nu_{n,k}(A_t)\ls t^{-kn}$. On the other hand from Corollary
\ref{cor:measure-ball-grass} we have $\nu_{n,k}(B_\rho(E,
\varepsilon))\gr (\varepsilon/c_1)^{k(n-k)}$. Choosing $t>1$ such
that $t^{-kn}= (\varepsilon/c_1)^{k(n-k)}$, that is $t\simeq
1/\varepsilon^{1-\frac{k}{n}}$, we conclude that $A_t^c \cap B_\rho(E,
\varepsilon)\neq \emptyset$ and the result follows. \prend

\medskip

\noindent {\it Proof of Theorem \ref{thm:main-result}.} It
follows from Theorem \ref{thm:main-2} applied for $\rho=d$ and
Theorem \ref{KM}. Moreover, if $L_{\mu} \ls C$ by \eqref{eq:3.13} we can take $k\ls n-1$.  \prend

\medskip

Note that the proof shows that Theorem \ref{thm:main-result} holds true if we replace
the distance $d$ with any other distance $\rho$ as in Theorem \ref{thm:Grass-entropy}.  
Let $k=\lambda n$ for some
$\lambda\in (0,1)$ and choose $\beta\simeq L_\mu$.
We have also proved that for any $E\in G_{n,\lambda n}$ and every $\varepsilon\in
(0,1)$ there exists $F\in G_{n,\lambda n}$ with $d(E,F)<\varepsilon$
and \begin{equation}
L_{\pi_F\mu}<\frac{CL_\mu}{\varepsilon^{1-\lambda}}.\end{equation}
One should compare the above inequality with the following (optimal) pointwise estimate: For
any isotropic log-concave measure $\mu$ on $\mathbb R^n$, every
$\lambda\in (0,1)$ and every $F\in G_{n,\lambda n}$ one has
\begin{equation}\label{eq:3.15} L_{\pi_F\mu}< (CL_\mu)^{1/\lambda}.
\end{equation} To see this recall the fact that for any isotropic
log-concave probability measure $\mu$ on $\mathbb R^n$ there exists an isotropic
convex body $T$ in $\mathbb R^n$ with the properties $L_T\simeq
L_\mu $ and
\begin{align*} \frac{L_{\pi_F\mu}}{L_\mu}\simeq |T\cap
F^\perp|^{1/k},
\end{align*} for all $F\in G_{n,k}$ (see \cite[Proposition 2.1, Lemma 5.8]{DP} for
details). Using the entropy estimate $N(L_TD_n, T)\ls (C_1L_T)^n$
and the Rogers-Shephard inequality from \cite{RS}
\begin{align}
  |P_FT|  |T\cap F^\perp|\ls {n\choose k} |T|
\end{align} for every $F\in G_{n,k}$
the estimate \eqref{eq:3.15} easily follows. In order to see that
\eqref{eq:3.15} is optimal consider an isotropic probability measure $\mu$ on
$\mathbb R^{\lambda n}$ and the measure $\nu$ on $\mathbb R^n$ with
$\nu=\mu\otimes \gamma_{(1-\lambda)n}$, where $\gamma_m$ is the
standard Gaussian on $\mathbb R^m$. Then, one can check that
$(e^{-1}L_\nu)^{1/\lambda} \ls  L_\mu=L_{\pi_{\mathbb R^{\lambda n}}
\nu}$.

We summarize the above discussion in the following:
\begin{proposition}
 \noindent Let $\mu$ be an isotropic $\log$-concave probability measure in $\mathbb R^{n}$ and $\lambda \in (0,1)$ and let  $E\in G_{n,\lambda n}$. Then
 \begin{equation}
  L_{\pi_E\mu}< (c_1L_{\mu})^{\frac{1}{\lambda}}
 \end{equation} and the inequality is sharp up to the constant $c_1$.
However, for every $\varepsilon >0$, there exists $F\in G_{n,\lambda n} $ such that $d(E,F) \ls\varepsilon$ and
 \begin{equation}
  L_{\pi_F\mu}<\frac{c_{2}L_\mu}{\varepsilon^{1-\lambda}},
 \end{equation}
 where $c_{1}, c_{2}>0$ are absolute constants.

\end{proposition}

\section{On the dependence on $\varepsilon$ in Theorem \ref{thm:main-result}}

In this section we discuss the dependence on $\varepsilon$ in Theorem \ref{thm:main-2} and we conclude with some remarks on the quantity $q_{v}(\mu, \beta)$.

We show that an improvement on the dependence on this parameters would imply the Hyperplane conjecture.

\noindent {\bf \S1.}  As we mentioned on the introduction, any improvement to
$o(1/\varepsilon)$ on the dependence on $\varepsilon$ in Theorem
\ref{thm:main-result}, would imply an affirmative answer to the
Hyperplane conjecture. More precisely we have the following:

\begin{quote} {\bf Assumption.} {\it There exist $\alpha \in (0,1)$ and positive
integers $k_n<n$, $k_n\to \infty$ with the following property: For
all $n$, for all isotropic log-concave probability measures $\mu$ on $\mathbb
R^n$, for all $\varepsilon \in (0,1)$ and for every $E\in G_{n,k_n}$
there exists $F\in G_{n,k_n}$ with $d(E,F)<\varepsilon$ and
$L_{\pi_F\mu}< C/\varepsilon^{\alpha}$.}
\end{quote}

\noindent Then we  prove that:

\begin{proposition}\label{prop:sol-under-assump} With the above assumption for
any $n$, any isotropic measure $\mu$ on $\mathbb R^n$ satisfies
$L_\mu <C$ for some absolute constant $C>0$.
\end{proposition}

\noindent In this paragraph, in view of \eqref{eq:3.2}, we define
the isotropic constant as: $$L_\nu:=|Z_m(\nu)|^{-1/m},$$ for any
isotropic log-concave probability measure $\nu$ on $\mathbb R^m$.

\medskip

For the proof of Proposition \ref{prop:sol-under-assump} we shall
need some lemmas. We start with the next stability result of the
isotropic constant of marginals with respect to the distance $d$.

\begin{lemma}
  \label{lem:local-isomor-via-O(n)} Let $K$ be a centrally symmetric convex body in $\mathbb
  R^n$ with $t=d_G(K, B_2^n)$. Let $E, F\in G_{n,k}$ with $d(E, F)=d$. Then, there exists
  $U\in O(n)$ such that $U(E)=F$ and
  \begin{equation}
    (1+td)^{-1}P_FK\subseteq U(P_EK)\subseteq (1+td) P_FK.
  \end{equation} In particular,
  \begin{equation}
    |P_EK|^{1/k}\ls (1+td)|P_FK|^{1/k}.
  \end{equation}
\end{lemma}

\noindent {\it Proof.} We consider $U\in O(n)$ such that
$d=\|I-U\|$. Let $\theta \in S_F$. Then, $U^\ast \theta=\phi\in S_E$
therefore we have $\|\theta-\phi\|_2\ls d$. We may write:
\begin{align*}
  \frac{h_{P_FK}(\theta)}{h_{U(P_EK)}(\theta)}=\frac{h_K(\theta)}{h_K(\phi)}\ls
  1+\frac{h_K(\theta-\phi)}{h_K(\phi)} \ls 1+\frac{dR(K)}{r(P_EK)},
\end{align*} where $R(\cdot), r(\cdot)$ are the circumradius and
inradius respectively. Similarly, we have that:
\begin{align*}
\frac{h_{U(P_EK)}(\theta)}{h_{P_FK}(\theta)}=\frac{h_K(\phi)}{h_K(\theta)}\ls
1+\frac{h_K(\theta -\phi)}{h_K(\theta)}\ls 1+\frac{dR(K)}{r(P_FK)}.
\end{align*} Since $r(P_FK), r(P_EK)\gr r(K)$ and $d_G(K,B_2^n)=R(K)/r(K)$ the result follows. \prend

\medskip

Applying Lemma \ref{lem:local-isomor-via-O(n)} for $K=Z_k(\nu)$ and
using the modified definition of the isotropic constant we arrive at
the following:

\begin{proposition}\label{prop:stability-iso}
Let $\nu$ be an isotropic $log$-concave probability measure in $\mathbb R^n$ and let $1\ls
k\ls n-1$. For any $E, F\in G_{n,k}$ we have:
  \begin{equation}
     L_{\pi_E\nu}/ L_{\pi_F\nu} \ls 1+ d_G(Z_k(\nu), B_2^n) d(E,F).
  \end{equation}
\end{proposition}

\smallskip

\noindent The above estimate also implies that the {\it length of the
gradient} (see \cite[Chapter 3]{Led} for a definition) of the
function $h:(G_{n,k},d)\to \mathbb R$ with $F\stackrel{h}\longmapsto
\log L_{\pi_F\nu}$ at $F$ is bounded by $d_G(Z_k(\nu), B_2^n)$, that
is $|\nabla h|(F)\ls d_G(Z_k(\nu),B_2^n)$ or
\begin{align}   |\nabla \log L_{\pi_F\nu}|\ls d_G(Z_k(\nu), B_2^n).
\end{align} In particular, the function $F\mapsto \log L_{\pi_F\nu}$
is $d_G(Z_k(\nu),B_2^n)$-Lipschitz with respect to $d$.

\smallskip

The next step is for a given measure $\mu$ on $\mathbb R^k$
($k=k_n$) to construct a measure $\nu$ on $\mathbb R^n$ such that
the geometric distance of $Z_k(\nu)$ to $B_2^n$ to be at most
$L_\mu$. We write $L_m=\sup_\nu L_\nu$ where the superemum is taken
over all isotropic $\log$-concave probability measures on $\mathbb R^m$.

\begin{proposition}\label{prop:exist-special-meas}
  For any $1\ls k<n$ there exists an isotropic $log$-concave probability measure $\mu_{1}$ on $\mathbb
  R^k$ with $L_{\mu_{1}}\simeq L_k$ and
  an isotropic $\log$-concave  measure $\mu_{2}$ on $\mathbb R^n$ such that $\pi_{\mathbb R^{k}}\mu_{2}= \mu_{1}$ and
\begin{align}\label{eq:4.7}
  \frac{c_1}{L_k}\sqrt{k}B_2^n \subseteq Z_k(\mu_{2})
\subseteq c_2\sqrt{k}B_2^n.
\end{align}
\end{proposition}

The next two lemmas will be needed for the proof of Proposition
\ref{prop:exist-special-meas}.

\begin{lemma}
 \label{lem:convol-cartesian-Z-k} The $L_q$-centroid bodies enjoy the following properties:
 \begin{itemize}
  \item [\rm (i)] If $\nu_1,\nu_2$ are probability measures in $\mathbb R^k$ and at least one of them is
  symmetric, then for all $q\gr 1$ we have:
  \begin{align}
   Z_q(\nu_1\ast
    \nu_2) \subseteq Z_q(\nu_1)+Z_q(\nu_2) \subseteq 2Z_q(\nu_1\ast
    \nu_2).
  \end{align}

  \item [\rm (ii)] If $\mu,\nu$ are probability measures on $\mathbb R^k$ and $\mathbb R^m$ respectively and at
  least one of them  is symmetric, then for all $q\gr 1$:
  \begin{align}
   Z_q(\mu\otimes \nu) \subseteq Z_q(\mu)\times Z_q(\nu) \subseteq 2 Z_q(\mu\otimes \nu) .
  \end{align}

 \end{itemize}

\end{lemma}

\noindent {\it Proof (Sketch).} We prove the second statement (see also Proposition 6.2 in \cite{KM}), the
first can be derived similarly - see \cite[Lemma 3.3]{GPV}. Since
for any $(x,y)\in \mathbb R^k\times \mathbb R^m$ we may write:
\begin{align*}
 h_{Z_q(\mu\otimes \nu)}(x,y) =\left(\int_{\mathbb R^k}\int_{\mathbb R^m} |\langle x,z_1\rangle+ \langle y,z_2\rangle|^q \, d\mu(z_1)\, d\nu(z_2)\right)^{1/q},
\end{align*}
the left-hand side inclusion follows from Minkowski's inequality
applied on the corresponding product space $(\mathbb R^k\times
\mathbb R^n, \mu \otimes \nu)$ for the functions $u(z_1,z_2)=\langle
x, z_1\rangle$ and $v(z_1,z_2)=\langle y, z_2\rangle$ and the fact
that for any two convex bodies $K,L$ we have $h_{K\times
L}(x,y)=h_K(x)+h_L(y)$. For the right-hand side inclusion we use the
symmetry to write:
\begin{align*}
 h_{Z_q(\mu\otimes \nu)}(x,y)=
 \left(\int_{\mathbb R^k}\int_{\mathbb R^m} \frac{|\langle x,z_1\rangle+ \langle y,z_2\rangle|^q+|\langle x,z_1\rangle- \langle y,z_2\rangle|^q}{2}
 \, d\mu(z_1)\, d\nu(z_2)\right)^{1/q}.
\end{align*}
Applying the elementary inequality:
$$|u+v|^q+|u-v|^q\gr |u|^q+|v|^q,$$ for all $u,v \in \mathbb R, \; q\gr 1$
we obtain: \begin{align*} h_{Z_q(\mu\otimes \nu)}(x,y) \gr
\left(\frac{h_{Z_q(\mu)}^q(x)+h_{Z_q(\nu)}^q(y)}{2}\right)^{1/q}.\end{align*}
The concavity of $t\mapsto t^{1/q}$ completes the proof.\prend

\medskip

Given an isotropic probability measure $\mu$ on $\mathbb R^k$ and any $\xi\in
(0,1)$ we define the measure $\mu_\xi$ with density function:
\begin{align}
  f_{\mu_\xi}(x):=\int_{\mathbb R^k} f_{\mu}(\sqrt{1-\xi^2}x-\xi
  y)g_k(\xi x+\sqrt{1-\xi^2}y)\, dy, \quad x\in \mathbb R^k,
\end{align} where $g_k(z)=(2\pi)^{-k/2}e^{-\|z\|_2^2/2}$ is the
density of the standard $k$-dimensional Gaussian measure $\gamma_k$. We need to adapt the definition of $M$-position for convex
bodies of V. Milman (see \cite{Mil}) in the setting of probability measures. We say that an isotropic $\log$-concave measure in $\mathbb R^{n}$ is in
{\it $M$-position with constant $A>0$}  if
the body $K:=L_{\mu}Z_k(\mu)$ satisfies: $|K+D_k|^{1/k}\ls A$.

\smallskip

Next lemma describes some properties of the measure
$\mu_\xi$.

\begin{lemma}
  \label{lem:Gauss-convol} Let $\mu$ be an isotropic $\log$-concave probability measure on $\mathbb R^k$.

  \begin{itemize}
  \item [1.] For any $\xi \in (0,1)$ the measure $\mu_\xi$ is
  log-concave and isotropic on $\mathbb R^k$ and \begin{equation}
  L_{\mu_\xi}\lesssim
  \min\left\{\frac{L_{\mu}}{\sqrt{1-\xi^2}},\frac{1}{\xi}
  \right\}.\end{equation}

  \item [2.] If $\mu$ is in $M$-position with constant $A>0$
  then we have:
  \begin{align}
    L_{\mu_\xi}\gtrsim \frac{1}{A}\min\left\{\frac{L_{\mu}}{\sqrt{1-\xi^2}},\frac{1}{\xi}
  \right\}.
  \end{align}
\item [3.] If $Z_k(\mu)\subseteq D_k$ then for any $\xi\in (0,1)$ we have:
  \begin{align}
   c_2 \xi D_k \subseteq Z_k(\mu_\xi)\subseteq c_3 D_k.
  \end{align}
  \end{itemize}
\end{lemma}

\noindent {\it Proof.} The log-concavity follows from
Pr\'{e}kopa-Leindler inequality \cite{Prek}. The isotropicity is
straightforward and follows from the fact that $\mu$ and $\gamma_k$
are isotropic. We may write:
\begin{align*}
  L_{\mu_\xi}\simeq f_{\mu_\xi}(0)^{1/k} &= \left(\int_{\mathbb R^k}
  f_{\mu}(-\xi y) g_k(\sqrt{1-\xi^2}y)\,dy \right)^{1/k} \\
  &\ls \|f_\mu\|_\infty^{1/k} \left(\int_{\mathbb R^k}
  g_k(\sqrt{1-\xi^2}y)\,
  dy\right)^{1/k}=\frac{L_{\mu}}{\sqrt{1-\xi^2}}.
\end{align*} Arguing similarly for $\gamma_k$, we can conclude that $L_{\mu_\xi} \ls C_1 \min\{L_\mu/\sqrt{1-\xi^2},
1/\xi\}$. For the inverse estimate we employ the information that
$\mu $ is in $M$-position. Considering the case $\xi\ls
\sqrt{1-\xi^2}/L_{\mu}$ we may write:
\begin{align*}
  |\sqrt{1-\xi^2}Z_k(\mu)+\xi
  D_k|^{1/k}\ls
  \frac{\sqrt{1-\xi^2}}{L_{\mu}}|K+D_k|^{1/k}\ls A
  \frac{\sqrt{1-\xi^2}}{L_{\mu}}.
\end{align*} In the case where $\xi\gr
\sqrt{1-\xi^2}/L_{\mu}$ we also have that
\begin{align*}
|\sqrt{1-\xi^2}Z_k(\mu)+\xi
  D_k|^{1/k}\ls A \max\{\sqrt{1-\xi^2}/L_\mu , \xi\}.
\end{align*} Therefore, using the fact that $Z_k(\mu_\xi)\simeq
\sqrt{1-\xi^2}Z_k(\mu)+\xi Z_k(\gamma_k)$ from Lemma
\ref{lem:convol-cartesian-Z-k} we obtain:
\begin{align*}
  L_{\mu_\xi}^{-1}\simeq |Z_k(\mu_\xi)|^{1/k}&\simeq |\sqrt{1-\xi^2}Z_k(\mu)+\xi
  Z_k(\gamma_k)|^{1/k} \lesssim A
  \max\left\{\frac{\sqrt{1-\xi^2}}{L_{\mu}},\xi\right\}
\end{align*} where we have also used the fact that $Z_k(\gamma_k)\simeq
D_k$. The last assertion follows again from Lemma \ref{lem:convol-cartesian-Z-k}(i). \prend

\medskip

\noindent {\it Proof of Proposition \ref{prop:exist-special-meas}.}
Let $1\ls k<n$ and let $\mu$ be an isotropic probability measure on $\mathbb
R^k$ with maximal isotropic constant. One can build (see \cite{DP}
for a construction) a new isotropic measure $\mu_0$ in $\mathbb R^k$
such that $L_{\mu_0}\simeq L_\mu$, $Z_k(\mu_0)\subseteq c_1 D_k$ and
$\mu_0$ is in $M$-position with absolute constant $A>0$. Selecting
$\xi\simeq L_\mu^{-1}$ and considering the measure $\mu_\xi$ induced by $\mu_0$ and
Gaussian we readily see that $L_{\mu_\xi}\simeq
L_\mu\simeq L_k$ and
\begin{align}\label{eq:4.12}
\frac{c_4}{L_\mu}\sqrt{k}B_2^k \subseteq Z_k(\mu_\xi) \subseteq
c_5\sqrt{k}B_2^k ,
\end{align} from Lemma \ref{lem:Gauss-convol}. Finally, consider the probability
measure $\mu_{2}:=\mu_\xi\otimes \gamma_{n-k}$. Clearly $\pi_{\mathbb R^{k}}\mu_{2}= \mu_\xi$. From Lemma
\ref{lem:convol-cartesian-Z-k} we know that
\begin{align}\label{eq:4.13}
  Z_k(\mu_{2}) \simeq Z_k(\mu_\xi)\times Z_k(\gamma_{n-k}).
\end{align} Since $Z_k(\gamma_{n-k})\simeq \sqrt{k}B_2^{n-k}$ the
inclusions in \eqref{eq:4.7} follow directly if we combine
\eqref{eq:4.12} and \eqref{eq:4.13} by setting $\mu_{1}=\mu_{\xi}$. \prend

\medskip

\noindent {\it Proof of Proposition \ref{prop:sol-under-assump}.}
Let $n\gr 1$ and $k=k_n<n$ the corresponding integer from the
assumption. Consider the measure $\mu_{1}$ on $\mathbb R^k$ and $\mu_{2} \in \mathbb R^{n}$ as given by Proposition
\ref{prop:exist-special-meas}. Then for $E=\mathbb R^k\times \{\mathbf
0\}$ and for $\varepsilon \simeq L_{k}^{-1}$ the assumption yields a
subspace $F\in G_{n,k}$ with $L_{\pi_F\mu_{2}}\ls
C\varepsilon^{-\alpha}$. By Proposition \ref{prop:stability-iso} we
obtain:
\begin{align}
  L_k\simeq L_{\mu_{1}}=L_{\pi_E\mu_{2}}\ls (1+c\varepsilon L_k)L_{\pi_F\mu_{2}}\ls C'
  L_{k}^{\alpha}
\end{align} and the claim is proved.  \prend

\medskip

\medskip

\noindent {\bf \S2.} The auxiliary parameter $q_v$ is one of the
many parameters have been introduced so far for the study of
isotropic log-concave measures. In \cite{Pa1} the parameter
$q_\ast(\mu)$ was introduced for proving sharp large deviation
estimates for the Euclidean norm with respect to a log-concave
measure. In \cite{Pa2} the parameter $q_{-c}(\mu,\delta), \,
\delta>1$ was introduced for the study of small ball probability
estimates. In \cite{KM} and \cite{Vri} local (hereditary) version of
these parameters was introduced for a unified approach to the
Hyperplane conjecture. We will not provide here all the definitions. The results of \cite[Theorem 1.2]{KM} and \cite[Theorem 1.1]{Vri} show that the quantity $q_{v}$ is larger than
the hereditary parameters. Moreover it is not hard one to construct examples (assuming that the Hyperplane conjecture is false) that the $q_{v}$ parameter is much larger than
the hereditary ones. This is not the case when one compares with the $q_{-c}$ parameter. For $\delta > 1$, $q_{-c}(\mu,\delta)$ is defined as the
largest $p$ such that $I_{-p}(\mu) \gr I_{2}(\mu)/\delta$, where $I_{q}(\mu):= \left(\int_{ \mathbb R^{n}} \|x\|_{2}^{q} d \mu(x) \right)^{\frac{1}{q}}$, $-(n-1)\ls q <\infty$.
It is known (see \cite{Pa2} Proposition 4.6) that for $k\ls q_{-c}(\mu,\delta)$
\begin{align*}
\left(\int_{G_{n,k}} L_{\pi_F\mu}^k \, d\nu_{n,k}(F)\right)^{\frac{1}{k}} \ls c\delta.
\end{align*}
So \eqref{eq:3.13}  implies that $q_{-c}(\mu,
C\beta)\gtrsim q_v(\mu, \beta)$ for any isotropic log-concave probability
measure $\mu$. If one could prove that for all $n$ and for all
log-concave, isotropic measures $\mu$ in $\mathbb R^n$ we have
$q_v(\mu,\beta)\gtrsim q_{-c}(\mu, \beta)$ for all $\beta>1$ the
Hyperplane conjecture would follow: if $\mu$ is an isotropic log-concave measure
on $\mathbb R^n$ then we can build the isotropic log-concave measure
$\nu=\mu\otimes \gamma_m$ in $\mathbb R^{n+m}$ where $\gamma_m$ is
the standard Gaussian and $m\simeq n/\log L_\mu$. Note that $L_\nu
\simeq L_\mu^{\frac{n}{n+m}}L_{\gamma_m}^{\frac{m}{n+m}}\gr
c_1L_\mu$. Moreover, we have $I_{-k}(\nu)\gr I_{-k}(\gamma_m)$ for
all $1\ls k\ls m-1$ which shows that $q_{-c}(\nu, \sqrt{\log
L_\mu})\gr c_2 n/\log L_\mu$. Then, by definition of $q_v$ we may
write:
$$L_\nu\ls c_3\sqrt{\log L_\mu}\sqrt{\frac{n+m}{q_v(\nu,
\sqrt{\log L_\mu})}}\ls c_4\sqrt{\log
L_\mu}\sqrt{\frac{n}{q_{-c}(\nu,\sqrt{\log L_\mu})}}\ls c_5{\log
L_\mu}.$$
Moreover the quantities $q_{\ast}$ and $q_{-c}$ are equivalent (up to constants) (see Proposition 2.5 in \cite{Pa4}) for ``truncated" isotropic
measures (supported on ball of radius of order $\sqrt{n}$). So the quantity $q_{v}(\mu, \beta)$ is on general larger than
the ``hereditary" quantities, smaller than the quantity $q_{-c}(\mu, \beta)$ and if one would prove that $q_{-c}$ and $q_{-v}$ are comparable the Hyperplane conjecture would follow.


\bigskip

\vspace{.5cm}
\noindent \begin{minipage}[l]{\linewidth}
  Grigoris Paouris: {\tt grigoris@math.tamu.edu}\\
  Petros Valettas: {\tt petvalet@math.tamu.edu}\\
  Department of Mathematics, Mailstop 3368\\
  Texas A\&M University\\
  College Station, TX 77843-3368\\

\end{minipage}

\end{document}